\pdfoutput=1
\RequirePackage{ifpdf}
\ifpdf 
\documentclass[pdftex]{sigma}
\else
\documentclass{sigma}
\fi

\usepackage{mathtools}
\usepackage{tikz-cd}
\usetikzlibrary{decorations.pathmorphing}
\usepackage{mathrsfs}
\usepackage{microtype}

\DeclareMathOperator{\Z}{\mathbb{Z}}
\DeclareMathOperator{\R}{\mathbb{R}}
\DeclareMathOperator{\Gr}{Gr}
\DeclareMathOperator{\SO}{SO}
\DeclareMathOperator{\Sym}{S}
\DeclareMathOperator{\tr}{tr}
\DeclareMathOperator{\I}{\mathcal{I}}
\DeclareMathOperator{\contact}{\mathcal{C}}

\renewcommand{\d}{{\rm d}}

\newcommand{\cinfinity}[1]{C^\infty\big(#1\big)}

\renewcommand{\mod}[1]{\quad\big(\mathrm{mod\enspace} #1\big)}
\newcommand{\pDeriv}[2]{\frac{\partial #1}{\partial #2}}

\newcommand{\w}{{\mathchoice{\,{\scriptstyle\wedge}\,}{{\scriptstyle\wedge}}
 {{\scriptscriptstyle\wedge}}{{\scriptscriptstyle\wedge}}}}

\newcommand{\M}[1]{M_0^{(#1)}}
\newcommand{\Ip}[1]{\I^{(#1)}}
\newcommand{\thetas}[1]{\theta_{#1}}
\newcommand{\thetanot}{\thetas{\varnothing,0}}
\newcommand{\pis}[1]{\pi_{#1}}
\newcommand{\omegas}[1]{\omega_{(#1)}}
\DeclareMathOperator{\pwf}{F}
\DeclareMathOperator{\swf}{\mathcal{F}}
\DeclareMathOperator{\wt}{\mathrm{wt}}
\DeclareMathOperator{\pwt}{\mathrm{pwt}}
\renewcommand{\dh}{\ensuremath{\hspace{1pt} \d_h}}
\newcommand{\dv}{\ensuremath{\hspace{1pt} \d_v}}

\newcommand{\barh}[1]{{}\mkern3mu\overline{\mkern-3mu H}^{#1}}
\newcommand{\cartTheta}{\mathcal{J}}
\newcommand{\cMinfty}{C^\infty(M)}
\newcommand{\jets}{J^2}
\newcommand{\Omegabar}{\overline{\Omega}}
\newcommand{\jet}[1]{J^{#1}\big(\R^{n+1},\R\big)}

\numberwithin{equation}{section}

\newtheorem{Theorem}{Theorem}[section]
\newtheorem{Corollary}[Theorem]{Corollary}
\newtheorem{Lemma}[Theorem]{Lemma}
\newtheorem{Proposition}[Theorem]{Proposition}
{ \theoremstyle{definition}
\newtheorem{Definition}[Theorem]{Definition}
\newtheorem{Example}[Theorem]{Example}
\newtheorem{Remark}[Theorem]{Remark} }

\begin{document}

\newcommand{\arXivNumber}{1810.02346}

\renewcommand{\PaperNumber}{047}

\FirstPageHeading

\ShortArticleName{Conservation Laws for a Class of Second-Order Parabolic Equations II}

\ArticleName{Geometry and Conservation Laws for a Class\\ of Second-Order Parabolic Equations II:\\ Conservation Laws}

\Author{Benjamin B. MCMILLAN}

\AuthorNameForHeading{B.B.~McMillan}

\Address{University of Adelaide, Adelaide, South Australia}
\Email{\href{mailto:benjamin.mcmillan@adelaide.edu.au}{benjamin.mcmillan@adelaide.edu.au}}

\ArticleDates{Received March 17, 2020, in final form April 27, 2021; Published online May 11, 2021}

\Abstract{I consider the existence and structure of conservation laws for the general class of~evolutionary scalar second-order differential equations with parabolic symbol. First I~calculate the linearized characteristic cohomology for such equations. This provides an~auxiliary differential equation satisfied by the conservation laws of a given parabolic equation. This is used to show that conservation laws for any evolutionary parabolic equation depend on at most second derivatives of solutions. As a corollary, it is shown that the only evolutionary parabolic equations with at least one non-trivial conservation law are of~Monge--Amp\`ere type.}


\Classification{35L65; 58A15; 35K10; 35K55; 35K96}

\section{Introduction}

The celebrated contributions of Noether began the systematic study of conservation laws for non-linear partial differential equations. More recently, Vinogradov introduced cohomological tools to the calculation of conservation laws, with the $\mathcal{C}$-spectral sequence of~\cite{Vinogradov:C-SpectralSequenceI} and~\cite{Vinogradov:C-SpectralSequenceII}. Much work has been developed from the $\mathcal{C}$-spectral sequence by Vinogradov and subsequent authors, but particularly relevant here is the work of Bryant and Griffiths in~\cite{Characteristic_Cohomology_I}, where they translate the $\mathcal{C}$-spectral sequence to the setting of exterior differential systems. Their approach is particularly amenable to considerations of the coordinate invariant geometry of a given class of differential equation.

In~a follow-up paper~\cite{Characteristic_Cohomology_II}, Bryant and Griffiths applied their technique to the class of non-linear second-order scalar parabolic equations for one unknown function of 2 variables. They defined the class of 7-dimensional exterior differential systems that are locally equivalent to such parabolic equations and showed that, in this case, conservation laws are in bijection with functions satisfying an auxiliary differential equation on the 7-manifold. Their 7-manifold is essentially the space of~2-jets of solutions, so an immediate consequence of their theorem is that any conservation law in this context depends on at most second derivatives of solutions. This should be contrasted with other classes of differential equation, of which there are many well known examples whose conservation laws depend on arbitrarily many derivatives of solutions~-- one famous example being the KdV equation. Bryant and Griffiths also showed that a 2-variable parabolic equation with at least one non-trivial conservation law is necessarily of Monge--Amp\`ere type.

There are new phenomena for parabolic equations in three or more variables. For example, it~is no longer the case that every parabolic equation can be put into evolutionary form. In~her thesis~\cite{Clelland_Thesis}, Clelland studied the exterior differential systems corresponding to evolutionary parabolic equations in 3 variables. She showed that in this case too, conservation laws depend on at most second derivatives of solutions and only Monge--Amp\`ere equations have non-trivial conservation laws.

In~this paper, I extend these results to evolutionary parabolic equations in any number of~variables. I show that for any evolutionary parabolic equation, conservation laws are functionals depending on at most second derivatives of solutions. As a corollary, if an evolutionary parabolic equation has at least one non-trivial conservation law, then it is necessarily Monge--Amp\`ere (Theorem~\ref{thm:parabolics with conservation laws are MA}).

The jump from 3 to all degrees of freedom relies on the author's recent progress on the geometry of parabolic differential equations.
In~\cite{McMillan:ParabolicsI} (henceforth, Part I), I introduce the general class of exterior differential systems that model second-order parabolic differential equations.
This puts the class of parabolic equations in a more geometric form, which I then study using Cartan's method of equivalence, introducing several local invariants of interest.

The local invariants derived from the equivalence problem~-- and their geometric interpretation~-- are crucial to the calculations in this paper.
Indeed, the extended Goursat invariants introduced in~Part~I allow for a simplification of the structure equations of evolutionary parabolic type exterior differential systems.
This reduction makes the calculations done here significantly more tractable.
The Monge--Amp\`ere invariants, also introduced in Part~I, provide a geometric characterization of Monge--Amp\`ere parabolic equations.
To wit, the existence of a non-trivial conservation law on an evolutionary parabolic equation allows one to show that a component of the Monge--Amp\`ere invariants vanishes identically, which is enough to show that the equation is Monge--Amp\`ere.

To be clear, this strategy is the same as was applied by the previous authors in dimen\-si\-ons~2 and~3.
The difficulty is that the complexity of the equivalence problem grows quickly; the number of local invariants is of order 4 in dimension.
This complexity can be mitigated in higher dimensions with a careful application of representation theory.
Indeed, the local invariants take values in various representations of~$\SO(n)$, so that we may consider the invariants in families.
From this perspective, the Goursat and the Monge--Amp\`ere invariants are treated the same in all dimensions, in much the same way as the Riemannian curvature tensor of Riemannian geometry~is.

To expand on this point, in any Cartan equivalence problem there is an associated Lie group that controls the local invariants.
The Lie group $ G $ associated to the geometry of parabolic equations is somewhat complicated to write down, but it contains a copy of~$\SO(n)$.
(Here $n+1$ is the number of independent variables in the given parabolic equation.)
The local invariants of a parabolic equation are equivariant functions from a $ G $-structure to a particular representation of~$ G $.
The local invariants can then be categorized according to which subrepresentations they take values in, and thus treated as unitary objects.
It is worth noting that in 2 and 3 dimensions the representation theory is obscured simply by low rank considerations, and there it is simpler to treat each scalar invariant separately.

\looseness=-1 The first several sections of this paper are more or less expository, devoted to recalling known results and tools, but with an emphasis on parabolic equations. In~Section~\ref{sec:background}, I~recall the geometric background on parabolic systems that will be applicable to conservation laws. The material in~that section is derived from results in~\cite{McMillan:ParabolicsI}. In~particular, I recall \emph{parabolic systems}, the exterior differential systems that are locally equivalent to parabolic equations.
Also of note are the geometric characterizations of evolutionary parabolic equations and of Monge--Amp\`ere parabolic equations.

In~Section~\ref{sec:prolongation of M_0}, I work out the structure equations for the infinite prolongation of a parabolic system $(M_0,\I_0)$. This is an infinite-dimensional replacement of~$M_0$ that has the same solutions as $M_0$, but sees information about arbitrarily many derivatives. This technical step is necessary to deal with the (a priori) possibility that conservation laws can depend on high derivatives of~solutions. To experts, there will be no surprises here, just a direct calculation depending on~the symbol type of a parabolic system.

Sections~\ref{sec:characteristic cohomology}, \ref{sec:weight filtrations}, and~\ref{sec:dv and dh} recall tools of~\cite{Characteristic_Cohomology_I}, but with an emphasis on parabolic equations.
Effort has been made throughout to present the material in such way that they may provide a useful guide for other examples.

In~Section~\ref{sec:conservation laws}, these tools are used to prove Theorem~\ref{thm:conservation laws linear part}, which describes the canonical form that any conservation law takes, and its dependence on a single function on the infinite prolongation of~$M_0$. This function satisfies an auxiliary differential equation on the infinite prolongation, so~that in principle, one needs to solve an infinite-dimensional PDE to find conservation laws.

In~Section~\ref{sec:Conservation laws of strong parabolic}, I prove
Theorem~\ref{thm:Jacobi potential lives on M_0}, which states that the defining function of a conservation law is in fact defined on $M_0$, instead of on the infinite prolongation. Since points of~$M_0$ are the 2-jets of solutions, this result can be restated as follows.
\begin{Theorem}
 Any conservation law for any evolutionary parabolic equation depends on at most second derivatives of solutions.
\end{Theorem}
This means in particular that the classification problem of finding all conservation laws for a~parabolic system is a finite-dimensional problem, in distinct contrast to other symbol classes.

Using the reduction of conservation laws to $M_0$ it is then fairly simple to prove Corollary~\ref{thm:parabolics with conservation laws are MA}, which can be paraphrased as follows.
\begin{Theorem}
 If an evolutionary parabolic equation has at least one non-trivial conservation law, then it is of Monge--Amp\`ere type in the neighborhood where the conservation law does not vanish.
\end{Theorem}

\section{Background}\label{sec:background}

I begin by briefly recalling some geometric results that will be needed here. For more details, and proofs, see~\cite{McMillan:ParabolicsI}.

In~coordinates, a scalar, $(n+1)$-variable, second-order \emph{weakly parabolic equation} is a differential equation for a single unknown function $u$ of~$n+1$ variables $x^0 , \dots , x^n$, the differential equation taking the form
\begin{gather}\label{eq: generic second order diffeq}
 F\bigg(x^a,u,\pDeriv{u}{x^a},\pDeriv{^2 u}{x^a\partial x^b}\bigg) = 0, \qquad a, b = 0,\dots , n ,
\end{gather}
subject to the constraint that the linearization at the 2-jet of any solution has parabolic symbol.

Here and throughout, ``space-time'' indices $a,b,\dots$ will range from $0$ to $n$, while ``spatial'' indices $i,j,\dots$ will range from $1$ to $n$. The Einstein summation convention will be used without further comment. In~fact, because all of the vector spaces that will be associated to a parabolic equation have a well defined \emph{spatial trace}, it will be convenient to apply the \emph{spatial Einstein summation} convention, so that, for example,
\begin{gather*}
\pDeriv{^2u}{x^i\partial x^i} := \sum_{i=1}^n \pDeriv{^2u}{x^i\partial x^i} .
\end{gather*}
The distinction between trace and spatial trace will be made by using (respectively) repeated space-time indices $a, b, \dots$ and spatial indices $i, j, \dots$.

A parabolic equation is \emph{evolutionary} if there is a choice of coordinates so that it is in~evolutionary form
\begin{gather*}
 \pDeriv{u}{x^0} = G\bigg(x^a,u,\pDeriv{u}{x^i},\pDeriv{^2 u}{x^i\partial x^j}\bigg) .
\end{gather*}
It should be noted that a generic weakly parabolic equation cannot be put into evolutionary form for any change of coordinates, even locally. Indeed, the extended Goursat invariants provide the obstruction to finding such coordinates. The Goursat invariants in turn provide the geometric condition that characterizes evolutionary parabolic equations. This condition is described below, and is crucial to the calculations to follow.

The following definition describes the exterior differential systems that are locally equivalent to scalar second-order weakly parabolic equations.
\begin{Definition}\label{def:weakly parabolic system}
 A \emph{weakly parabolic system} in $n+1$ variables is a
 $(2n+2+(n+1)(n+2)/2)$-dimensional\footnote{This is 1 less than the dimension of~$J^2(\R^{n+1},\R)$, corresponding to the fact that a parabolic equation is defined by a single equation on $J^2(\R^{n+1},\R)$.} exterior differential system $(M_0,\I_0)$ such that any point has a neighborhood equipped with a spanning set of 1-forms
 \begin{gather*}
 \thetas{\varnothing},\qquad \thetas{a},\qquad\omega^a,\qquad\pis{ab} = \pis{ba}
 \end{gather*}
 that satisfy:
 \begin{enumerate}\itemsep=0pt
 \item
 The forms $\thetas{\varnothing}$, $\thetas{a}$ generate $\I_0$ as a differential ideal.
 \item
 The structure equations
 \begin{gather*}
 \d\thetas{\varnothing} \equiv -\thetas{a}\w\omega^a \ \mod{\thetas{\varnothing}},
 \\
 \d\thetas{a} \equiv -\pis{ab}\w\omega^b \mod{\thetas{\varnothing},\thetas{b}} .
 \end{gather*}
 \item
 The parabolic symbol relation (see comment above about spatial Einstein summation)
 \begin{gather*}
 \pis{ii} \equiv 0 \mod{\thetas{\varnothing},\thetas{a},\omega^a} .
 \end{gather*}
 \end{enumerate}
\end{Definition}

\begin{Remark}
 There is a natural association from weakly parabolic equations to parabolic systems.
 This association is such that
 the (graphs of) solutions of the parabolic equation are in~bijection with the solution submanifolds of the corresponding parabolic system.
 For the reader's convenience, I briefly recall the essential facts here, but refer to
 Example~1.4 of~\cite{McMillan:ParabolicsI} for more details.

 The space of 2-jets $J = J^2\big(\!\R^{n+1},\R\!\big)$ is isomorphic as a manifold to $ \R^{n+1} \!\times\! \R \!\times\! \R^{n+1} \!\times\! \Sym\big(\!\R^{n+1}\!\big) $.
 A choice of coordinate $ x^a $ on $ \R^{n+1} $ and $ u $ on $ \R $ determines jet coordinates $ x^a$, $u$, $p_a$, $p_{ab}=p_{ba} $ on~$ J $, where the $p_a$ correspond to the first derivatives of~$u$ with respect to $x^a$ and $p_{ab}$ to the second derivatives.
 These coordinates may be used to define a~coframing of~$ J $,
 \begin{gather}\label{eq: generic second order coframing}
 \hat\theta_\varnothing := \d u - p_a\d x^a, \qquad
 \hat\theta_a := \d p_a - p_{ab}\d x^b, \qquad \d x^a, \qquad\d p_{ab} .
 \end{gather}
 The coframing here depends on the choice of coordinates, but the following structure equations hold for any such choice,
 \begin{align*}
 \d\hat\theta_\varnothing & \equiv - \d p_a \w \d x^a \mod{\hat\theta_{\varnothing}},
 \\
 \d\hat\theta_a & \equiv \d p_{ab} \w \d x^b \mod{\hat\theta_{\varnothing}, \hat\theta_{b}} ,
 \end{align*}
 and these structure equations are reflected in condition 2 of Definition~\ref{def:weakly parabolic system}.

 Now, a second order differential equation such as equation~\eqref{eq: generic second order diffeq} defines, in a clear way, a function $ F(x^a,u,p_a,p_{ab}) $ on $ J $, and assuming $ F $ is sufficiently non-degenerate, the zero set $ M_0 = F^{-1}(0) $ is a manifold of codimension~1.
 If $ \Sigma \hookrightarrow M_0 $ is a submanifold for which the ideal $ \I_0 = \big\{ \hat\theta_{\varnothing}, \hat\theta_{b} \big\} $ pulls back to zero and the pullback of~$ \d x^0 \w\cdots\w \d x^n $ is nowhere vanishing,
 then $ \Sigma $ is locally the 2-jet graph of a solution to the differential equation $F$.
 Such $ \Sigma $ are called \emph{solution submanifolds} of the exterior differential system $ (M_0, \I_0) $, and it is straightforward to see that the 2-jet graph of a solution to the differential equation is a solution submanifold.

 The coframing~\eqref{eq: generic second order coframing} pulls back to $ M_0 $, and the equation $ \d F = 0 $ determines a single relation between the forms of the coframing.
 In~particular, for $ F $ a weakly parabolic equation, there exists at every point of~$ M_0 $ a change of coframe from~\eqref{eq: generic second order coframing} to one satisfying the conditions of~Definition~\ref{def:weakly parabolic system}.

 Conversely, every small enough neighborhood in a parabolic system arises in this way from a~weakly parabolic differential equation.
 See~\cite[Theorem~5.10]{BCGGG} or the proof of Theorem~5.3 in~\cite{McMillan:ParabolicsI} for details.
 On the other hand, there are parabolic systems which don't have a global embedding into $\jets$, such as the parabolic system modelling mean curvature flow, taken up in Example~1.6 of~\cite{McMillan:ParabolicsI}.
\end{Remark}

\begin{Remark}\label{rmk: parabolic G structure}
 Every spanning set of 1-forms as in the definition is adapted to the parabolic symbol type of the system. As such, they are called \emph{parabolic $($extended$)$ coframings}. It is worth noting that they are not strictly coframings, because of the parabolic symbol relation. Nonethless, an extended parabolic coframing may be used to define a coframing: for any point $x \in M_0$, the extended coframing defines a linear injection
 \begin{gather*}
 T_xM \xrightarrow{\qquad} \R \oplus \R^{n+1} \oplus \big(\R^{n+1}\big)^\vee \oplus \Sym^2\R^{n+1} .
 \end{gather*}
 By construction, post-composing this with the quotient map annihilating the spatial trace element in $\Sym^2 \R^{n+1}$, we get an isomorphism of~$T_xM$ with a canonical vector space~-- a~legimate coframe. The issue then is that this quotient space does not have a geometrically natural choice of basis. This unnatural choice can be avoided by using extended parabolic coframings.

 Parabolic coframes can be used to define a $G$-structure in the standard manner, and it is in~this context that the geometric invariants are developed.
\end{Remark}

\begin{Remark}\label{rmk: Cartan system}
 It follows from the structure equations of any parabolic coframing that the 1-form~$\thetas{\varnothing}$ is uniquely defined up to rescaling by a function. As a consequence, the Cartan system of~$\thetas{\varnothing}$ is well defined independent of a choice of parabolic coframing. It can be directly computed from the structure equations that Cartan system of~$\thetas{\varnothing}$ is the Frobenius ideal
 \begin{gather*}
 \cartTheta = \{ \thetas{\varnothing}, \thetas{a}, \omega^a \} .
 \end{gather*}
 This ideal will help to define vertical derivatives of functions in Section~\ref{sec:dv and dh}.
\end{Remark}

Evolutionary equations are geometrically privileged in the class of all parabolic equations, and more can be said of them.
Indeed, there is geometric test for evolutionarity: it was shown in~Theorem 5.4 of~\cite{McMillan:ParabolicsI} that
a given parabolic equation can be put in evolutionary form if and only if its associated parabolic system has a choice of parabolic coframing for which
\begin{gather*}
\d \omega^0 \equiv 0 \mod{\omega^0} .
\end{gather*}
This in turn is equivalent to the refined structure equations
\begin{gather}\label{eq:refined evolutionary structure equations}
 \d\thetas{i} \equiv -\pis{ia}\w\omega^a \mod{ \thetas{\varnothing},\thetas{j} }
\end{gather}
for $i, j = 1 , \dots , n$.
Furthermore, because the sub-principal symbol of an evolutionary parabolic is non-vanishing, the coframing may be chosen so that
\begin{gather}\label{eq:evolutionary parabolic sub-principal symbol}
 \pis{ii} = \thetas{0} .
\end{gather}

In~fact, an evolutionary parabolic system has parabolic coframings satisfying the structure equations
\begin{gather*}
 \d\thetas{\varnothing} \equiv - \thetas{a}\w\omega^a \mod{\thetas{\varnothing} }, \\
 \d\thetas{i} \equiv - \pis{ia}\w\omega^a \mod{ \thetas{\varnothing}, \thetas{i} }, \\
 \d\thetas{0} \equiv - \pis{0a}\w\omega^a \mod{ \thetas{\varnothing}, \thetas{a} }, \\
 \d\omega^0 \equiv 0 \mod{\omega^0 }, \\
 \d\omega^i \equiv 0 \mod{\thetas{\varnothing}, \thetas{a}, \omega^a } .
\end{gather*}
Furthermore, from $\d^2 \thetas{i} = 0$, one readily computes that
\begin{gather*}
 \d\pis{ij} \equiv 0 \mod{ \thetas{\varnothing}, \thetas{a}, \omega^a, \pis{kl} },
 \\
 \d\pis{i0} \equiv 0 \mod{ \thetas{\varnothing}, \thetas{a}, \omega^a, \pis{kl}, \pis{k0} } .
\end{gather*}

Monge--Amp\`ere parabolic equations are even more geometrically privileged than evolutionary parabolics, and have evident deep connections with conservation laws. An evolutionary parabolic equation is \emph{Monge--Amp\`ere} if it is of the form
\begin{gather*}
 \pDeriv{u}{x^0} = \sum_{\substack{I, J \subseteq \{1,\dots,n\} \\ |I|=|J|}} A_{I,J}\bigg(x^a,u,\pDeriv{u}{x^i}\bigg) H_{I,J} ,
\end{gather*}
where the $A_{I,J}$ are arbitrary smooth functions and $H_{I,J}$ is the row $I$, column $J$ minor sub-determinant of the Hessian of~$u$. Monge--Amp\`ere equations are quantifiably less complex than the generic case: their corresponding exterior differential systems have a natural de-prolongation to an exterior differential system of dimension only $2n+3$.
\begin{Definition}
 A \emph{quasi parabolic Monge--Amp\`ere system} in $n+1$ variables is a $(2n+3)$-di\-men\-sio\-nal exterior differential system $(M_{-1},\I_{-1})$ such that $\I_{-1}$ is locally generated by a 1-form $\thetas{\varnothing}$ and an~$(n+1)$-form $\Upsilon$ satisfying:
 \begin{enumerate}\itemsep=0pt
 \item $\thetas{\varnothing}$ is maximally non-integrable:
 \begin{gather*}
 \thetas{\varnothing} \w (\d\thetas{\varnothing})^{n+1} \neq 0 .
 \end{gather*}
 \item At each point of~$M_{-1}$,
 \begin{gather*}
 \Upsilon \not\equiv 0 \mod{\thetas{\varnothing}, \d\thetas{\varnothing}} .
 \end{gather*}
 \item There is (locally) a 1-form $\omega^0$ on $M_{-1}$ independent of~$\thetas{\varnothing}$ and so that
 \begin{gather*}
 \omega^0\w\Upsilon \equiv 0 \mod{\thetas{\varnothing}, \d\thetas{\varnothing}} ,
 \end{gather*}
 and any other such 1-form is a linear combination of~$\thetas{\varnothing}$ and $\omega^0$.
 \end{enumerate}
\end{Definition}

In~\cite{McMillan:ParabolicsI}, I developed the local invariants that provide an effective test for a parabolic system to have a Monge--Amp\`ere deprolongation. Briefly, for an evolutionary parabolic system $M_0$ and any evolutionary parabolic coframing, there is a $\mathfrak{co}(n)$-valued 1-form $\big(\beta_i^j\big) = \big(\delta_i^j\beta_{\tr} + \mathring{\beta}_i^j\big)$ (so~$\mathring{\beta}_i^j = -\mathring{\beta}_j^i)$, a symmetric traceless matrix valued 1-form $\big(\xi_i^j\big)$, and a 1-form $\kappa_{\varnothing}$ so that
\begin{gather}
 \d\thetas{i} \equiv - \beta_i^j\w\thetas{j} - \xi^j_i\w\thetas{j} - \pis{ia}\w\omega^a \mod{\thetas{\varnothing},\Lambda^2 \I_0}, \nonumber
 \\
 \d\omega^i \equiv -\big(\delta^i_j \kappa_{\varnothing} - \beta^i_j\big)\w\omega^j + \xi^i_j\w\omega^j \mod{\thetas{\varnothing}, \thetas{a}, \omega^0} .\label{eq: evolutionary MA structure equations}
\end{gather}
There are functions $V^{jkl}_i$ on $M_0$ so that
\begin{gather*}
\xi_i^j \equiv V_i^{jkl}\pis{kl} \mod{\thetas{\varnothing}, \thetas{a}, \omega^a } .
\end{gather*}
The functions $V_i^{jkl}$ generically take values in the $\mathfrak{co}(n)$ representation $\Sym^2 \big(\Sym^2_0 \R^n\big)$, but a coframe reduction may be made absorbing the $ik$-trace, so that $V_i^{jkl}$ takes values in $\Sym^2_0 \big(\Sym^2_0 \R^n\big)$. These functions are the secondary\footnote{As the name suggests, there are primary Monge--Amp\`ere invariants. They vanish automatically for evolutionary parabolics.} Monge--Amp\`ere invariants, which determine whether~$M_0$ comes from a parabolic equation of Monge--Amp\`ere type. The following is part of Theorem~4.3 in~\cite{McMillan:ParabolicsI}.
\begin{Theorem}
 For a parabolic system $(M_0,\I_0)$ the following are equivalent:
 \begin{enumerate}\itemsep=0pt
 \item[$1.$] $M_0$ is locally EDS equivalent to a neighborhood in the prolongation of a quasi-parabolic Monge--Amp\`ere system.

 \item[$2.$] The Monge--Amp\`ere invariants $V_i^{jkl}$ of~$M_0$ take values in
 \begin{gather*}
 \mathfrak{b}_2 := \ker\big(\Sym^2_0 \big( \Sym^2_0 \R^n\big) \xrightarrow{\qquad} \Sym_0^3 \R^n \otimes \R^n\big).
 \end{gather*}
 \end{enumerate}
\end{Theorem}

It is worth noting that
\begin{gather*}
\Sym^2_0 \big( \Sym^2_0 \R^n\big) = \mathfrak{b}_2 \oplus \Sym_0^4 \R^n ,
\end{gather*}
so it suffices for the totally symmetric component of the secondary Monge--Amp\`ere invariant to vanish.
For the parabolic system $M_0$ associated to a parabolic equation, condition (1) is equivalent to the parabolic equation being of Monge--Amp\`ere type. On the other hand, once the conservation law theory has been worked out here, condition (2) will be immediately obvious for any parabolic system with at least one non-trivial conservation law.

\section[The prolongation of (M0, I0)]{The prolongation of~$\boldsymbol{(M_0,\I_0)}$}\label{sec:prolongation of M_0}

A priori, the conservation laws of a parabolic equation may depend on arbitrarily many derivatives of solutions. In~order to prove that this is not the case for evolutionary parabolic equations, it will be necessary to work on the infinite prolongation of~$(M_0, \I_0)$, whose definition I recall here.

Let $\jet{r}$ denote the manifold of~$r$-jets of functions from $\R^{n+1}$ to $\R$. A choice of linear coordinates $x^0,\dots, x^n$ on the domain $\R^{n+1}$ and coordinate $p_\varnothing$ on the codomain $\R$ induces natural coordinates on $\jet{r}$~-- for the symmetric multi-index $I$ in $\{0 , \dots , n\}$, the coordinate $p_I$ corresponds to the derivative $\frac{\partial^{|I|}p_\varnothing}{\partial x^I}$.

The jet manifolds may be arranged in a tower, where each $\jet{r+1}$ is a bundle over $\jet{r}$ with fiber isomorphic to the symmetric product $\Sym^{r+1}\big(\R^{n+1}\big)$. Furthermore, each jet space $\jet{r}$ has a natural exterior differential system structure, with the differential ideal~$\contact^{(r)}$ generated by the 1-forms\footnote{Here and throughout, $Ia$ will denote the symmetric multi-index obtained from $I$ by appending $a$.}
\begin{gather*}
\hat\theta_I = \d p_I - p_{Ia}\d x^a, \qquad |I| < r .
\end{gather*}

Let $(\jet{\infty}, \contact^{(\infty)})$ denote the inverse limit of this tower. On $\jet{\infty}$ there are natural differential operators, the \emph{total derivatives}, defined by the formal sum
\begin{gather*}
D_a = \partial_{x^a} + \sum_I p_{Ia}\partial_{p_I} .
\end{gather*}
These can be used to describe the operation of prolongation. For example, consider a second-order differential equation
\begin{gather*}
F\bigg(x^a, p_\varnothing, \pDeriv{p_\varnothing}{x^a}, \pDeriv{p_\varnothing}{x^a \partial x^b}\bigg) = 0
\end{gather*}
for the unknown scalar function $p_\varnothing$ of~$n+1$ variables, which defines a function
\begin{gather*}
F(x^a, p_\varnothing, p_a, p_{ab})
\end{gather*}
on $\jet{\infty}$ that clearly factors through $\jet{2}$. The first prolongation of this equation is, by definition, the system of~$n+2$ functions
\begin{gather*}
F^{(1)} \colon\ \jet{3} \xrightarrow{\qquad} \R\times\R^{n+1} ,
\end{gather*}
where $F^{(1)} = F\times(D_a F)$.
This process may be repeated inductively to define the prolonged systems $F^{(2)}, F^{(3)}, \dots$ on to $F^{(\infty)}$. It is important to note that solutions of~$F$ are in bijection with those of~$F^{(r)}$ for each $r$, so no information is lost at any step.

 Now consider a parabolic system $(M_0, \I_0)$, as defined above. Locally on $M_0$ there is a parabolic second-order equation $F$ and an embedding so that the diagram of EDS maps
 \begin{gather*}
 (M_0, \I_0) \xhookrightarrow{\quad \iota \quad} \big(\jet{2}, \contact^{(2)}\big) \xrightarrow{\quad F \quad} (\R, \{0\})
 \end{gather*}
 is exact, in the sense that $M_0 = F^{-1}(0)$ and $\I_0 = \iota^*\contact^{(2)}$.

The exterior differential system prolongations of~$M_0$ (see~\cite{BCGGG, Characteristic_Cohomology_I} for the intrinsic definition and more details) can be seen to fit into the diagram
\[
{\begin{tikzcd}
\vdots \arrow{d} & \vdots \ar[d] & \vdots \arrow[d]
\\
\big(M_0^{(2)}, \I_0^{(2)}\big) \arrow{d}\arrow[hookrightarrow]{r} & \big(\jet{4}, \contact^{(4)}\big) \ar[d]\ar[r, "F^{(2)}"] & \big(\R^{1 + n+1 + \binom{n+1}{2}}, \{0\}\big) \arrow[d]
\\
\big(M_0^{(1)}, \I_0^{(1)}\big) \ar[d]\arrow[hookrightarrow]{r} & \big(\jet{3}, \contact^{(3)}\big) \ar[d]\ar[r, "F^{(1)}"] & \big(\R^{1 + n+1}, \{0\}\big) \arrow[d]
\\
(M_0, \I_0) \arrow[hookrightarrow]{r} & \big(\jet{2}, \contact^{(2)}\big) \ar[r, "F"] & (\R, \{0\})
\end{tikzcd}}
\]
It is a key point that the left column may be calculated globally, without reference to local embeddings into jet space. In~particular, infinite prolongations exist even for exterior differential systems that have no global embedding into jet bundles, such as the mean curvature flow.
Furthermore, any $(n+1)$-dimensional solution submanifold in $M_0$ has a unique lift to $\M{r}$ for each $r$, so the solution manifolds of~$M_0$ are in bijection with those of~$\M{r}$ for each $r$.

Let $(M,\I)$ be the infinite prolongation of~$M_0$, the inverse limit of the prolongation tower. More precisely, the manifold $M$ is given by the inverse limit of underlying manifolds,
\begin{gather*}
M = \varprojlim \M{r} ,
\end{gather*}
and the ideal $\I$ is given by
\begin{gather*}
\I = \bigcup_{r=0}^\infty \Ip{r} .
\end{gather*}

The manifold $M$ is infinite-dimensional, but this will not cause any technical difficulty in~prac\-tice, because any conservation law factors through some finite level. In~particular, we are concerned with conservation laws of \emph{finite type}, which are represented by certain differential forms pulled back from a finite prolongation of~$M_0$.
Accordingly, it will suffice to consider finite-type functions on $M$, those that can be expressed as the pullback of a function on $\M{r}$ for some $r$.
The space of finite-type smooth functions is given by
\begin{gather*}
\cMinfty = \bigcup_{r=0}^\infty \cinfinity{\M{r}}
\end{gather*}
and the finite-type differential forms by
\begin{gather*}
\Omega^*(M) = \bigcup_{r=0}^\infty \Omega^*\big(\M{r}\big) .
\end{gather*}

It is worth noting that a function $A \in \cinfinity{\M{r}}$ is one which, when restricted to the $r$-jet graph of a solution $u(x)$ to the differential equation $F$, results in a functional of~$x$, of~$u$, and of~the derivatives of~$u$ up to order $r+2$. For example, a point of~$M_0$ is the 2-jet of a solution, and a function on $M_0$ is a functional depending on at most second derivatives of solutions.

Now for some linear algebra preliminaries, which will be used to describe the structure of~$(M,\I)$. Fix a vector space $W$ and subspace $W'$ of respective dimensions $n+1$ and $n$, and basis $e_a$ on $W$ so that
\begin{gather*}
W' = \R\{e_1 , \dots , e_n\} \subset W = \R\{e_0 , \dots , e_n\} .
\end{gather*}
Let
\begin{gather*}
W'' = W/W' \cong \R\{e_0\} .
\end{gather*}
This decomposition of~$W$ induces one on the symmetric powers of~$W$, so that
\begin{gather*}
\Sym^{r} W \cong \bigoplus_{s+t= r} \Sym^s( W') \otimes \Sym^t( W'') .
\end{gather*}
For notation, elements will be indexed with regards to this splitting, so that $\Sym^r(W)$ has basis given by the elements
\begin{gather*}
e_{I,t} := e_{i_1}{\circ}\cdots{\circ}e_{i_s}{\circ}(e_0)^{t}
\end{gather*}
for $I = (i_1\dots i_s)$ as $I$ ranges over symmetric multi-indices of length $|I|$ up to $r$, and $t = r -|I|$.

It will be convenient to relabel the adapted coframing $\omega^a$, $\thetas{\varnothing}$, $\thetas{a}$, $\pis{ab}$ of~$M_0$
so that it is consistent with the basis of~$\Sym^2 W$:
\begin{gather*}
 \thetas{\varnothing} \xmapsto{\quad} \thetas{\varnothing,0}, \qquad \pis{ij} \xmapsto{\quad} \pis{ij,0},
 \\
 \thetas{i} \xmapsto{\quad} \thetas{i,0},\qquad \ \,\pis{i0} \xmapsto{\quad} \pis{i,1},
 \\
 \thetas{0} \xmapsto{\quad} \thetas{\varnothing,1},\qquad \pis{00} \xmapsto{\quad} \pis{\varnothing, 2}.
\end{gather*}

The total symmetric product $\Sym_\bullet := \Sym_\bullet(W)$ of~$W$ is an algebra, with multiplication given by~the action of~$W$:
\begin{gather*}
 e_i e_{I,t} = e_{Ii,t} ,
 \\
 e_0 e_{I,t} = e_{Ii,t+1} .
\end{gather*}
Note that this is nothing but a free polynomial algebra with a special indeterminate picked out.
There is also an action of~$W^\vee$ on $\Sym_\bullet$, which is essentially the directional derivative. In~the given basis,
\begin{gather*}
 e^i e_{I,t} := \#(i,I) e_{I\backslash i,t}, \qquad
 e^0 e_{I,t+1} := (t+1) e_{I,t} ,
\end{gather*}
where $\#(i,I)$ is the number of times $i$ appears in $I$.
These can be used to define the \emph{spatial trace} operator on $\Sym_\bullet$, given by
\begin{gather*}
\tr = \sum_{i=1}^n e^ie^i .
\end{gather*}

Now, a parabolic system $(M_0, \I_0)$ is a linear Pfaffian system, and the general theory of such (see~\cite[Chapter~IV]{BCGGG}) determines the associated tableaux
\begin{gather*}
K \cong \Sym^2_0 (W') \oplus (W'\otimes W'') \oplus \Sym^2 (W'') \subset W \otimes W ,
\end{gather*}
where $\Sym^2_0 (W')$ is the traceless component of~$\Sym^2 (W')$.
The $r$-fold prolongation of~$K$ is by definition
\begin{gather*}
K^{(r)} = (K\otimes \Sym^r W) \cap \big(W \otimes \Sym^{r+1} W\big) ,
\end{gather*}
which is readily computed to be
\begin{gather*}
K^{(r)} \cong \bigoplus_{s = 0}^{2 + r} \Sym^s_0 (W') \otimes \Sym^{2 + r - s} (W'') .
\end{gather*}
Each $K^{(r)}$ is naturally identified with the kernel of the spatial trace
\begin{gather*}
 \tr \colon\ \Sym^{r+2}(W) \xrightarrow{\quad e^{i}e^i \quad} \Sym^{r}(W) .
\end{gather*}

This linear algebra determines the structure equations in the following proposition, whose proof is standard.
\begin{Proposition}[principal structure equations]\label{thm:principal structure equations}
 Near any point of the infinite prolonga\-tion~$(M, \I)$ of a parabolic system, there is a spanning set of $1$-forms
 \begin{gather*}
 \omega^a, \qquad \thetas{I,t} ,
 \end{gather*}
 where $I$ ranges over symmetric multi-indices and $t$ over non-negative integers, so that
 \begin{enumerate}\itemsep=0pt
 \item[$1.$] $\I = \{\thetas{I,t}\}$,

 \item[$2.$] $\I^{(r)} = \{\thetas{I,t} \colon |I| + t \le r + 1\}$,

 \item[$3.$] for any $\thetas{I,t}$,
 the principal structure equations
 \begin{gather}\label{eq:principal structure equations}
 \d\thetas{I,t} \equiv - \thetas{Ii,t}\w\omega^i - \thetas{I,t+1}\w\omega^0 \mod{\thetas{J,s} \colon |J| + s \le |I| + t}
 \end{gather}
 hold,

 \item[$4.$] for each $r\ge 0$, the pullback to $M$ of~$\Omega^1\big(\M{r}\big)$ is locally spanned by $\omega^a$ and the 1-forms of~$\I^{(r+1)}$, and

 \item[$5.$] for each fixed $I$ and $t$, the $ 1 $-forms in $\I$ are subject to the relation
 \begin{gather*}
 \sum_{i=1}^n \thetas{Iii,t} \equiv 0 \mod{\thetas{J,s} \mbox{ for which } |J| + s < |I| + t + 2} .
 \end{gather*}
 \end{enumerate}
\end{Proposition}

\begin{Definition}
 Any choice of~$\omega^a$, $\thetas{I,t}$ as in Proposition~\ref{thm:principal structure equations} is called a~\emph{parabolic coframing} of~$M$.
\end{Definition}

\begin{Remark}\label{rmk: infinitely prolonged parabolic G structure}
 As in Remark~\ref{rmk: parabolic G structure}, a parabolic coframing only defines a coframing of~$M$ after composing with the appropriate quotient. Nonetheless, the bundle of all parabolic coframes can be used to define a $G^{(\infty)}$-structure, where $G^{(\infty)}$ is the infinite group prolongation of~$G$. Essentially by construction, $G^{(\infty)}$ acts stably on the filtration of~$\Omega^*(M)$ by the ideals $\I^{(r)}$.
\end{Remark}

The theorem has the following corollary.
\begin{Corollary}\label{cor: factoring through}
 A function $A \in \cMinfty$ factors through $\cinfinity{\M{r}}$ if and only if
 \begin{gather*}
 \d A \equiv 0 \mod{\omega^a, \I^{(r+1)}} .
 \end{gather*}
\end{Corollary}
\begin{proof}
 The function $A$ is finite type, so factors through $\cinfinity{\M{R}}$ for sufficiently large $R$. Since the fibers of~$\M{R}$ are connected, $A$ factors through $\cinfinity{\M{r}}$ if and only if $\d A \in \Omega^1\big(\M{r}\big)$, which by the theorem holds if and only if
 \begin{gather*}
 \d A \equiv 0 \mod{\omega^a, \I^{(r+1)}} .
 \tag*{\qed}
\end{gather*}
\renewcommand{\qed}{}
\end{proof}

\section{The characteristic cohomology}\label{sec:characteristic cohomology}

With $(M, \I)$ the infinitely prolonged parabolic system described in the last section, let $\Omegabar^*$ be chain complex defined by the exact sequence
\begin{gather*}
0 \xrightarrow{\quad} \I^* \xrightarrow{\quad} \Omega^*(M) \xrightarrow{\quad} \Omegabar^* \xrightarrow{\quad} 0 .
\end{gather*}
Because $\I$ is a differential ideal, the differential of~$\Omega^*(M)$ descends to a differential $\dh$ on the quotient $\Omegabar^*$, which is therefore a graded commutative differential algebra.
This cdga is intimately connected with the behavior of solutions to $M$. The next definition follows Bryant and Griffiths~\cite{Characteristic_Cohomology_I}.
\begin{Definition}
 The \emph{characteristic cohomology} $\barh{q}$ of an infinitely prolonged parabolic sys\-tem~$(M,\I)$ is the cohomology of the complex $\Omegabar^*$,
 \begin{gather*}
 \barh{q} = H^q\big(\Omegabar^*,\dh\big) .
 \end{gather*}
\end{Definition}
Due to the symbol of~$(M,\I)$, it follows from the results in~\cite{Characteristic_Cohomology_I} that $\barh{n}$ and $\barh{n+1}$ are the only nontrivial characteristic cohomology groups. Furthermore,
$\barh{n}$ is naturally identified as the space of conservation laws of~$(M,\I)$:
\begin{Definition}
 The \emph{space of conservation laws} for a parabolic system is given by the degree $n$ characteristic cohomology,
 \begin{gather*}
 \overline{\mathscr{C}} = \barh{n}(M) .
 \end{gather*}
\end{Definition}

The cohomology of the differential operator $\dh$ is difficult to compute directly. However, the characteristic cohomology fits into a spectral sequence. In~turn, the first page of this spectral sequence can be computed using a second spectral sequence. This second spectral sequence is quite amenable to calculations, as the differential of it first page is linear over functions, allowing a pointwise computation determined entirely by the symbol of~$(M,\I)$.

To define the first spectral sequence, consider the filtration of~$\Omega^*(M)$ defined recursively by
\begin{gather*}
 \pwf^0 = \Omega^*(M), \\
 \pwf^{p+1} = \I\w \pwf^p .
\end{gather*}
Each level of the filtration is graded, with graded components
\begin{gather*}
\pwf^{p,q} = \pwf^p \cap \Omega^{p+q}(M) .
\end{gather*}
Then the bi-graded associated graded spaces are given by
\begin{gather*}
Gr^{p,q} := \pwf^{p,q}/\pwf^{p+1,q-1}.
\end{gather*}
$\I$ is formally Frobenius, so, for each $p\ge 0$, the exterior derivative descends to a well defined operator $\dh$. These complexes comprise page 0 of the filtration spectral sequence, so that
\begin{gather*}
E^{p,q}_0 = Gr^{p,q} .
\end{gather*}
It follows immediately from the definition that
\begin{gather*}
E^{0,q}_0 = \Omegabar^q ,
\end{gather*}
and thus
\begin{gather*}
E^{0,q}_1 = \barh{q} .
\end{gather*}

On the other hand,
\begin{gather*}
E_\infty \Rightarrow H^*(M) ,
\end{gather*}
so for any contractible neighborhood in $M$ the spectral sequence converges to the homology of~a~point. As we are concerned here with local conservation laws, assume henceforth that $M$ is contractible, restricting attention to a neighborhood if necessary.

It follows from the two-line theorem of Vinogradov~\cite{Vinogradov}, that for parabolic equations,
\begin{gather*}
E^{p,q}_1 = 0 \text{\quad for\quad } q < n ,
\end{gather*}
so the $E_1$ page is
\begin{equation*}
 {\begin{tikzcd}
 0 \ar[r] & {\barh{n+1}} \ar[r, "\dv"] & {E^{1,n+1}_1} \ar[r, "\dv"] & {E^{2,n+1}_1} \ar[r, "\dv"] & \cdots ,
 \\[-3ex]
 0 \ar[r] & {\barh{n}} \ar[r, "\dv"] & {E^{1,n}_1} \ar[r, "\dv"] & {E^{2,n}_1} \ar[r, "\dv"] & \cdots ,
 \\[-3ex]
 0 \ar[r] & 0 \ar[r] & 0 \ar[r] & 0 \ar[r] & \cdots.
 \end{tikzcd}}
\end{equation*}
Note that as an immediate consequence, the bottom row is exact, so that $\barh{n}$ is isomorphic to the kernel of~$\dv$ in $E^{1,n}_1$. This motivates the following definition, from~\cite{Characteristic_Cohomology_I}.
\begin{Definition}
 The \emph{space of differentiated conservation laws} is
 \begin{gather*}
 \mathscr{C} = \ker\big(E^{1,n}_1 \xrightarrow{\ \dv\ } E^{2,n}_1\big) .
 \end{gather*}
\end{Definition}
The operator $\dv$ provides an isomorphism between $\overline{\mathscr{C}}$ and $\mathscr{C}$, and the latter space may be computed in two steps: First one computes $E^{1,n}_1$, and then one computes the kernel of~$\dv$. It is the second spectral sequence that helps to compute $E_1^{1,n}$, which is explained in the next section.

\section{The weight filtrations}\label{sec:weight filtrations}

Bryant and Griffiths introduced a second filtration on $E_1^{p,q}$ when $p\ge 1$, the \emph{principal weight filtration}, that linearizes the calculation of each $E^{p,q}_1$.

\begin{Definition}
 A \emph{weight function} is a function $\wt\colon \Omega^*(M) \to \Z$ satisfying the following properties:
 \begin{enumerate}\itemsep=0pt
 \item $\wt(f) = 0$ for $f \in \cMinfty$.
 \item $\wt(\alpha\w\beta) = \wt(\alpha) + \wt(\beta)$ for $\alpha,\beta \in \Omega^*(M)$.
 \item $\wt(\alpha + \beta) = \max(\wt(\alpha),\wt(\beta))$.
 \end{enumerate}
\end{Definition}

\begin{Example}[principal weight filtration]
 Fix any parabolic coframing of~$M$ and consider the weight function $\pwt$ on $\Omega^*(M)$ uniquely specified by
 \begin{gather*}
 \pwt(\omega^a) = -1,
 \\
 \pwt(\thetas{I,t}) = |I| + t .
 \end{gather*}
 The weight function $\pwt$ descends to $\Gr^{p,*}$ for each $p>0$. As argued in~\cite{Characteristic_Cohomology_I}, $\pwt$ on $\Gr^{p,*}$ is independent of the choice of parabolic coframing. See also Remark~\ref{rmk: infinitely prolonged parabolic G structure}.

 For each integer $k$ and each $p \ge 1$, define
 \begin{gather*}
 \pwf^{p}_k =\big \{ \alpha \in \pwf^p \colon \pwt(\alpha) \le k\big \}/\pwf^{p+1} ,
 \end{gather*}
 the subspace in $\Gr^{p,*}$ of forms with principal weight less than $k$.
 The grading by form degree descends to $\pwf^p_k$, with graded components
 \begin{gather*}
 \pwf^{p,q}_k = \big\{ \alpha \in \pwf^{p,q} \colon \pwt(\alpha) \le k\big \}/\pwf^{p+1} .
 \end{gather*}
 The \emph{principal weight filtration} on $\Gr^{p,*}$ is defined by the sequence
 \begin{gather*}
 \cdots \subset \pwf^{p}_k \subset \pwf^{p}_{k+1} \subset \cdots \subset \Gr^{p,*} .
 \end{gather*}
 It follows immediately from equations~\eqref{eq:principal structure equations} that this filtration is stable under $\dh$.

 For each fixed $p>0$ there is a spectral sequence $\overline{E}$ associated to this filtration, with $\overline{E}_0$ page given by
 \begin{gather*}
 \overline{E}_0^{q,k} = \pwf^{p,q}_k/\pwf^{p,q}_{k-1} .
 \end{gather*}
 This spectral sequence converges to $E^{p,*}_1$.

 In~calculations, it will be convenient to abuse notation slightly, using $\pwf^p_k$ again to denote the preimage of~$\pwf^p_k$ under the mapping $\pwf^p \to \Gr^{p,*}$, that is,
 \begin{gather*}
 \pwf^p_k = \big\{ \alpha \in \pwf^p \colon \pwt(\alpha) \le k \big\} + \pwf^{p+1}
 \end{gather*}
 and
 \begin{gather*}
 \pwf^{p,q}_k = \big\{ \alpha \in \pwf^{p,q} \colon \pwt(\alpha) \le k \big\} + \pwf^{p+1} .
 \end{gather*}
 Then for any element $[\Phi] \in \overline{E}_0^{q,k}$ represented by a form $\Phi$, the choice of~$\Phi$ is uniquely defined modulo $\pwf^p_{k-1}$. For example, if $\dh\Phi$ is any representative of~$\dh[\Phi] \in \overline{E}_0^{q+1,k}$, then{\samepage
 \begin{gather*}
 \dh \Phi \equiv \d\Phi \mod{\pwf^p_{k-1}} .
 \end{gather*}
 The distinct usages occur in different contexts, and should not cause confusion.
 }

 Any element of~$\overline{E}_0^{q,k}$ can be represented by a linear combination of~$(p+q)$-forms
 \begin{gather*}
 \thetas{I_1,t_1}\w\cdots\w\thetas{I_p,t_p}\w\omega^{a_1}\w\cdots\w\omega^{a_{q}}
 \end{gather*}
 of principal weight exactly $k$ -- explicitly, those for which
 \begin{gather*}
 |I_1| + t_1 +\cdots + |I_p| + t_p - q = k .
 \end{gather*}

 The exterior derivative $\dh$ is $\cMinfty$-linear on $\overline{E}_0^{*,k}$, and treats the forms $\omega^a$ as constants. Indeed, by $\R$-linearity of~$\dh$, it suffices to check that
 \begin{gather*}
 \dh (f\thetas{I_1,t_1}\w\cdots\w\thetas{I_p,t_p}\w\omega^{a_1}\w\cdots\w\omega^{a_q})
 \\ \qquad
 {} \equiv f\dh(\thetas{I_1,t_1}\w\cdots\w\thetas{I_p,t_p})\w\omega^{a_1}\w\cdots\w\omega^{a_q} \mod{\pwf^{p}_{k-1}} ,
 \end{gather*}
 but this follows from the observation that $\dh$ strictly decreases weight for functions and for each~$\omega^a$.

 Corollary~\ref{cor: factoring through} can be restated in terms of the principal weight filtration:
 \begin{Corollary}\label{cor:factoring through by weight}
 A function $A \in \cMinfty$ factors through $\cinfinity{\M{r}}$ if and only if
 \begin{gather*}
 \d A \equiv 0 \mod{\omega^a, \pwf^1_{r+2}} .
 \end{gather*}
 \end{Corollary}
\end{Example}

\begin{Example}[sub-principal weight filtration]
 For parabolic systems, the principal weight filtration doesn't see ``lower order'' information, such as the sub-principal symbol. So, it will also be necessary to introduce the \emph{sub-principal weight filtration}, which sees further into the structure of~$(M,\I)$. The proof of Theorem~\ref{thm:Jacobi potential lives on M_0} relies heavily on the sub-principal weight filtration.

 Unlike the previous case, the sub-principal weight filtration on $\Gr^{p,*}$ will necessarily depend on a choice of parabolic coframing. However, it will momentarily be shown that there are more refined choices of parabolic coframing and that all such coframings define the same sub-principal weight filtration.

 Relative to a parabolic coframing, the \emph{sub-principal weight function} is the unique weight function $\wt$ such that
 \begin{gather*}
 \wt\big(\omega^i\big) = -1, \\
 \wt\big(\omega^0\big) = -2, \\
 \wt(\thetas{I,t}) = |I| + 2t .
 \end{gather*}

 For each integer $k$ and each $p \ge 1$, define the subspace
 \begin{gather*}
 \swf^p_k = \big\{ \alpha \in \pwf^p \colon \wt(\alpha) \le k \big\}/ \pwf^{p+1}
 \end{gather*}
 of~$\Gr^{p,*}$.
 The \emph{sub-principal weight filtration} on $\Gr^{p,*}$ is defined by the sequence
 \begin{gather*}
 \cdots \subset \swf^p_k \subset \swf^p_{k+1} \subset \cdots \subset \Gr^{p,*} .
 \end{gather*}
 As with the principal weight filtration, it will be useful to have the notation
 \begin{gather*}
 \swf^p_k = \big\{ \alpha \in \pwf^p \colon \wt(\alpha) \le k \big\} + \pwf^{p+1}
 \end{gather*}
 and
 \begin{gather*}
 \swf^{p,q}_k = \big\{ \alpha \in \pwf^{p,q} \colon \wt(\alpha) \le k \big\} + \pwf^{p+1} .
 \end{gather*}

 For a generic choice of parabolic coframing, the induced filtration is not automatically $\dh$-stable. So make the following definition.
 \begin{Definition}
 A parabolic coframing $\omega^a$, $\thetas{I,t}$ is a \emph{refined parabolic coframing} of~$M$ if it satisfies the refined structure equations
 \begin{gather}\label{eq:refined structure eqs}
 \d\thetas{I,t} \equiv -\thetas{Ii,t}\w\omega^i -\thetas{I,t+1}\w\omega^0 \mod{\swf^1_{|I|+2t-1}}.
 \end{gather}
 \end{Definition}

 Before turning to the proof that refined parabolic coframings exist for evolutionary parabolic systems, I state several useful properties of such coframings. Only the proof of the second property relies on the parabolic system being evolutionary.
 \begin{Proposition}
 Suppose given a refined parabolic coframing $\omega^a$, $\thetas{I,t}$ on an evolutionary parabolic system and corresponding sub-principal weight filtration $\swf$.
 \begin{enumerate}\itemsep=0pt
 \item[$1.$] The ideal $\swf^p_{S}$ is $\dh$-stable for each $p>0$ and all $S$.

 \item[$2.$] On $\swf^p_S / \swf^p_{S-1}$, the operator $\dh$ is $\cMinfty$-linear and treats the $\omega^a$ as constants.

 \item[$3.$] The higher symbol relations
 \begin{gather*}
 \thetas{Iii,t} \equiv \thetas{I,t+1} \mod{\swf^1_{I+2+2t-1}}
 \end{gather*}
 hold for all $I$ and $t$.

 \item[$4.$] Any other choice of refined parabolic coframing gives the same sub-principal weight filtration.
 \end{enumerate}
 \end{Proposition}
 \begin{proof}
 1. It suffices to check this for elements in $\swf^p_{S}$ of the form
 \begin{gather*}
 f \thetas{I_1,t_1}\w\cdots\w\thetas{I_p,t_p}\w\omega^{a_1}\w\cdots\w\omega^{a_{q}} .
 \end{gather*}
 The claim is clear from the structure equations~\eqref{eq:refined structure eqs} and from the fact that $\dh$ does not increase weight for functions and the $\omega^a$.

 2. This follows because $\dh$ strictly decreases sub-principal weight of functions and the $\omega^a$. In~particular, because the system is assumed evolutionary,
 \begin{gather*}
 \dh\omega^0 \equiv 0 \mod{\omega^0, \pwf^1} .
 \end{gather*}

 3. Equation~\eqref{eq:evolutionary parabolic sub-principal symbol}, pulled back to $M$, gives
 \begin{gather*}
 \thetas{ii,0} = \thetas{\varnothing,1} ,
 \end{gather*}
 providing the base case for induction.
 If
 \begin{gather*}
 \thetas{Iii,t} \equiv \thetas{I,t+1} \mod{\swf^1_{|I|+2+2t-1}}
 \end{gather*}
 holds, then applying $\d$ to both sides,
 \begin{gather*}
 - \thetas{Ijii,t}\w\omega^j - \thetas{Iii,t+1}\w\omega^0 \equiv -\thetas{Ij,t+1}\w\omega^j - \thetas{I,t+2}\w\omega^0 \mod{\swf^1_{|I|+2+2t-1}} ,
 \end{gather*}
 so an application of Cartan's lemma gives
 \begin{gather*}
 \thetas{Ijii,t} - \thetas{Ij,t+1} \equiv 0 \mod{\omega^a, \swf^1_{|I| + 2 + 2t}}, \\
 \thetas{Iii,t+1} - \thetas{I,t+2} \equiv 0 \mod{\omega^a, \swf^1_{|I| + 3 + 2t}} .
 \end{gather*}
 But the left hand side is contained in $\I$, so cannot have nonzero $\omega^a$ terms, and thus
 \begin{gather*}
 \thetas{Ijii,t} - \thetas{Ij,t+1} \equiv 0 \mod{\swf^1_{|I| + 2 + 2t}}, \\
 \thetas{Iii,t+1} - \thetas{I,t+2} \equiv 0 \mod{\swf^1_{|I| + 3 + 2t}} .
 \end{gather*}

4. Let $\tilde\omega^a$, $\tilde\theta_{I,t}$ be a second choice of refined parabolic coframing, with associated sub-principal weight filtration $\tilde\swf$.
 It suffices to show that for each weight $N$ and sub-principal weight $S$,
 \begin{gather*}
 \pwf^{1,0}_N \cap \swf^{1,0}_S = \pwf^{1,0}_N \cap \tilde\swf^{1,0}_S .
 \end{gather*}

For $N = 0$, there is nothing to show, as $\pwf^{1,0}_0$ is spanned by $\thetas{\varnothing,0}$ which is a multiple of~$\tilde\theta_{\varnothing,0}$.

 For each weight $N > 0$,
 \begin{gather*}
 \pwf^{1,0}_N \subset \swf^{1,0}_{2N} ,
 \end{gather*}
 for if $\thetas{I,t}$ is such that $|I|+t \le N$, then $|I|+2t \le 2N$. Likewise,
 \begin{gather*}
 \pwf^{1,0}_N \subset \tilde\swf^{1,0}_{2N} .
 \end{gather*}
 Thus, for each $N$,
 \begin{gather*}
 \pwf^{1,0}_N \cap \swf^{1,0}_{2N} = \pwf^{1,0}_N = \pwf^{1,0}_N \cap \tilde\swf^{1,0}_{2N} .
 \end{gather*}

 So, suppose $N$ is the first weight for which the claim fails, and then $S < 2N$ the largest sub-principal weight such that
 \begin{gather*}
 \pwf^{1,0}_N \cap \swf^{1,0}_S \neq \pwf^{1,0}_N \cap \,\tilde\swf^{1,0}_S .
 \end{gather*}
 Then there are some $\tilde\theta_{Ii,t} \in \pwf^{1,0}_N\cap\,\tilde\swf^{1,0}_S$ and functions $A_{Ii,t}^{J,s}$ for which
 \begin{gather*}
 \tilde\theta_{Ii,t} \equiv \sum_{\substack{|J|+s \le N \\ |J|+2s > S}} A_{Ii,t}^{J,s}\thetas{J,s} \not\equiv 0 \mod{\pwf^{1,0}_N\cap\swf^{1,0}_{S}}
 \end{gather*}
 and thus for $\tilde\theta_{I,t} \in \pwf^{1,0}_{N-1}\cap\tilde\swf^{1,0}_{S-1} = \pwf^{1,0}_{N-1}\cap\swf^{1,0}_{S-1}$,
 \begin{gather*}
 \d\tilde\theta_{I,t} \equiv -\tilde\theta_{Ii,t}\w\omega^i \equiv \sum A_{Ii,t}^{J,s}\thetas{J,s}\w\omega^i \mod{\omega^0, \pwf^{1,1}_{N-2}\cap\swf^{1,1}_{S-2}} .
 \end{gather*}
 The right hand side is not in $\pwf^{1,1}_{N-1}\cap\swf^{1,1}_{S-1}$. However, $\pwf^{1}_{N-1}\cap\swf^{1}_{S-1}$ is $\dh$-stable, so $\d\tilde\theta_{I,t}$ is in~$\pwf^{1,1}_{N-1}\cap\swf^{1,1}_{S-1}$. This contradiction proves the statement.
 \end{proof}

 \begin{Proposition}
 Let $(M,\I)$ be the infinite prolongation of an evolutionary parabolic system. There exist refined parabolic coframings.
\end{Proposition}
 \begin{proof}
 It will suffice to inductively prove the following statement for each $N > 0$: there is a~parabolic coframing such that for each $\thetas{I,t}$ with $\pwt(\thetas{I,t}) \le N$,
 \begin{gather}\label{eq:strongest structure eqs}
 \d\thetas{I,t} \equiv -\thetas{Ii,t}\w\omega^i -\thetas{I,t+1}\w\omega^0 \mod{\pwf^1_{|I|+t-1}\cap\swf^1_{|I|+2t-1}} .
 \end{gather}

 The base case $N = 1$ follows immediately from the pullback to $M$ of the refined structure equations~\eqref{eq:refined evolutionary structure equations} for an evolutionary parabolic system.

 Suppose the induction has been carried out up to $N-1$.
 Then, for each sub-principal weight~$k$, one has
 \begin{gather*}
 \d\big(\pwf^{1,1}_{N - 2}\cap\swf^{1,1}_{k}\big) \subset \pwf^{1,2}_{N - 2}\cap\swf^{1,2}_{k} .
 \end{gather*}
 This follows because, modulo $\pwf^2$, any element of~$\pwf^{1,1}_{N - 2}\cap\swf^{1,1}_{k}$ is a linear combination of terms $\thetas{I,t}\w\omega^a$ such that $\pwt(\thetas{I,t}) \le N - 1$.

 It follows from the principal structure equations of Proposition~\ref{thm:principal structure equations} and from
 \begin{gather*}
 \pwf^{1,1}_{N-1} = \big\{\dots, \thetas{\varnothing,N}\w\omega^i\big\} \subset \swf^{1,1}_{2N-1}
 \end{gather*}
 that
 \begin{align*}
 \d\theta_{0,N} & \equiv -\thetas{i,N}\w\omega^i - \thetas{0,N+1}\w\omega^0 \mod{\pwf^1_{N-1}}, \\
 & \equiv -\thetas{i,N}\w\omega^i - \thetas{0,N+1}\w\omega^0 \mod{\pwf^1_{N-1}\cap\swf^1_{2N-1}}
 \end{align*}
 which provides the base case for a second induction.

 Suppose then that~\eqref{eq:strongest structure eqs} holds for all $\thetas{J,s}$ with
 \begin{gather*}
 \pwt(\thetas{J,s}) = N , \qquad \wt(\thetas{J,s}) > N + S .
 \end{gather*}
 Consider $\thetas{Ii,t}$ of weights
 \begin{gather*}
 \pwt(\thetas{Ii,t}) = N, \qquad \wt(\thetas{Ii,t}) = N + S .
 \end{gather*}
 From the principal structure equations, there are functions so that
 \begin{gather*}
 \d\thetas{Ii,t} \equiv -\thetas{Iij,t}\w\omega^j - \thetas{Ii,t+1}\w\omega^0 - \sum G^{(J,s)}_{(Ii,t),a}\thetas{J,s}\w\omega^a \mod{\pwf^1_{N-1}\cap\swf^1_{N+S-1}} ,
 \end{gather*}
 the sum over indices $J, s, a$ such that
 \begin{gather*}
 \pwt(\thetas{J,s}\w\omega^a) \le N-1 \quad \mbox{ and } \quad \wt(\thetas{J,s}\w\omega^a) \ge N + S .
 \end{gather*}
 From the first induction hypothesis,
 \begin{gather*}
 \d\thetas{I,t} \equiv - \thetas{Ii,t}\w\omega^i - \thetas{I,t+1}\w\omega^0 \mod{\pwf^{1,1}_{N-2}\cap\swf^{1,1}_{N+S-2}}
 \end{gather*}
 and thus (using now the second induction)
 \begin{gather*}
 0 = \d^2\thetas{I,t} \equiv \sum G^{(J,s)}_{(Ii,t),a}\thetas{J,s}\w\omega^a\w\omega^i \mod{\pwf^{1,2}_{N-2}\cap\swf^{1,2}_{N+S-2}} ,
 \end{gather*}
 so that
 \begin{gather*}
 G^{(J,s)}_{(Ii,t),j} = G^{(J,s)}_{(Ij,t),i} \qquad \text{and} \qquad G^{(J,s)}_{(Ii,t),0} = 0 .
 \end{gather*}

 From this it follows that one may make the coframe modification
 \begin{gather*}
 \thetas{Iij,t} \xmapsto{\qquad} \thetas{Iij,t} - G^{(J,s)}_{(Ii,t),j}\thetas{J,s} ,
 \end{gather*}
 which absorbs all of the $G^{J,s}_{Ii,t,j}$ while preserving the principal structure equations and the previous induction steps.
 \end{proof}

\end{Example}

\section{Horizontal and vertical derivatives}\label{sec:dv and dh}

The differentials $\dh$ and $\dv$ in the characteristic spectral sequence have natural geometric interpretations, which I describe here. They are essentially the same operators as defined in~\cite{Characteristic_Cohomology_I}.

By definition, the operator $\dh$ is the restriction of~$\d$ to the associated graded $\Gr^{p,*}$. More expli\-ci\-tly, given an equivalence class $[\alpha] \in \Gr^{p,*}$ represented by $\alpha \in \pwf^p$, the horizontal derivative~$\dh[\alpha]$ is represented by $\d\alpha \in \pwf^p$, so that
\begin{gather*}
\dh\alpha \equiv \d\alpha \mod{\pwf^{p+1}} .
\end{gather*}

For any element $[\varphi] \in \pwf^0/\pwf^1$, the operator $\dh$ is \emph{horizontal} with respect to any solution manifold $\Sigma$:
\begin{gather*}
(\dh \varphi)|_\Sigma = \d(\varphi|_\Sigma) .
\end{gather*}
In~particular, a function $A \in \cMinfty$ is a functional on solution manifolds $\Sigma$, from which we may define new functionals $D_a A$ so that
\begin{gather*}
\dh A \equiv (D_a A) \omega^a \mod{\pwf^1} .
\end{gather*}
Then from
\begin{gather*}
(\dh A)|_\Sigma = \d(A|_\Sigma)
\end{gather*}
it follows that the operators $D_a$ are locally the restriction to $M$ of the total derivatives $D_a$ defined on $J^\infty\big(\R^{n+1},\R\big)$.

Now consider the \emph{vertical derivative} $\dv$. Given a form $\alpha \in \pwf^{p}$ such that $\dh\alpha \equiv 0$, the deriva\-tive~$\d\alpha$ is an element of~$\pwf^{p+1}$. As such, define the vertical derivative of~$\alpha$ by
\begin{gather*}
\dv\alpha \equiv \d\alpha \mod{\pwf^{p+2}} .
\end{gather*}
Then, the horizontal derivative of this vanishes,
\begin{gather*}
\dh\dv\alpha \equiv \dh\d\alpha \equiv \d^2\alpha = 0 \mod{\pwf^{p+2}} .
\end{gather*}
This last observation will be useful, because it is often easier to compute $\dh\dv$ than $\dh$, as the operator $\dv$ aids one in filtering by weights.

It will be useful for calculations to extend the operator $\dv$ to all functions. This is accomplished in an invariant manner using the Cartan system $\cartTheta = \{\thetanot, \thetas{a}, \omega^a\} = \big\{\omega^a, \pwf^{1,0}_1\big\}$ defined in~Remark~\ref{rmk: Cartan system}, so that
\begin{gather*}
\dv A \equiv \d A \mod{\cartTheta} .
\end{gather*}
The following lemma shows how to compute the vertical derivative of the directional total derivatives of a function.
\begin{Lemma}\label{thm:(lem) vertical derivative of directional derivative, sub-principal weights}
 Let $M$ be the infinite prolongation of an evolutionary parabolic system. For any function $A \in \cMinfty$ with $w = \wt(\d A)$, there are functions $A_{I,s}$, not all zero, for which
 \begin{gather*}
 \d A \equiv (D_aA)\omega^a + \sum_{|I|+2s = w} A_{I,s}\thetas{I,s} \mod{\swf^1_{w-1}} .
 \end{gather*}
 Then
 \begin{gather*}
 \dv(D_i A) \equiv \sum A_{I,s}\thetas{Ii,s} \mod{\cartTheta, \pwf^1_2, \swf^1_{w}} ,
 \\
 \dv(D_0 A) \equiv \sum A_{I,s}\thetas{I,s+1} \mod{\cartTheta, \pwf^1_2, \swf^1_{w+1}} .
 \end{gather*}
 Furthermore, if $w \ge 4$, then
 \begin{gather*}
 \dv(D_i A) \equiv \sum A_{I,s}\thetas{Ii,s} \mod{\cartTheta, \swf^1_{w}} ,
 \\
 \dv(D_0 A) \equiv \sum A_{I,s}\thetas{I,s+1} \mod{\cartTheta, \swf^1_{w+1}} .
 \end{gather*}
\end{Lemma}
\begin{proof}
 From the structure equations for $\omega^i$ and $\omega^0$,
 \begin{gather*}
 \d\omega^a \equiv 0 \mod{\Lambda^2\cartTheta, \pwf^1_1} .
 \end{gather*}
 Then
 \begin{align*}
 0 = \d^2 A \equiv{}& \d \bigg((D_a A)\omega^a + \sum_{|I|+2s = w} A_{I,s}\thetas{I,s}\bigg) \mod{\swf^1_{w-1}}, \\
 \equiv{}& \dv (D_a A)\omega^a - \sum_{|I|+2s = w} A_{I,s}\big(\thetas{Ii,s}\omega^i + \thetas{I,s+1}\omega^0\big) \mod{\Lambda^2 \cartTheta, \pwf^1_1, \swf^1_{w-1}} ,
 \end{align*}
 so that the result follows from an application of Cartan's lemma.
 The last claim follows from the observation that $\pwf^{1,0}_2 \subseteq \swf^{1,0}_4$.
\end{proof}

The previous lemma holds with principal weight replacing sub-principal weight. The proof is, \textit{mutatis mutandis}, the same, so omitted.
\begin{Lemma}
 Let $M$ be the infinite prolongation of a parabolic system. For any function $A \in \cMinfty$, with $N = \pwt(\d A) \ge 2$, there are functions $A_{I,s}$, not all zero, for which
 \begin{gather*}
 \d A \equiv (D_aA)\omega^a + \sum_{|I|+s = N} A_{I,s}\thetas{I,s} \mod{\pwf^1_{N-1}} .
 \end{gather*}
 Then
 \begin{gather*}
 \dv(D_i A) \equiv \sum A_{I,s}\thetas{Ii,s} \mod{\cartTheta, \pwf^1_N} ,
 \\
 \dv(D_0 A) \equiv \sum A_{I,s}\thetas{I,s+1} \mod{\cartTheta, \pwf^1_N} .
 \end{gather*}
\end{Lemma}

\begin{Corollary}
 If $f \in \cinfinity{\M{r}}$, then $D_af \in \cinfinity{\M{r+1}}$.
\end{Corollary}
\begin{proof}
 For a function $f \in \cinfinity{\M{r}}$, one has by Corollary~\ref{cor:factoring through by weight},
 \begin{gather*}
 \dv f \equiv 0 \mod{\cartTheta, \pwf^1_{r+2}} ,
 \end{gather*}
 and thus by the previous lemma,
 \begin{gather*}
 \dv D_a f \equiv 0 \mod{\cartTheta, \pwf^1_{r+3}} ,
 \end{gather*}
 so that $D_a f \in \cinfinity{\M{r+1}}$.
\end{proof}

\section{Conservation laws of parabolic systems}\label{sec:conservation laws}
The first step in calculating the space of conservation laws for a parabolic system is to com\-pute~$E^{1,n}_1$. The statement of the following Theorem uses the omitted index notation, wherein for an~anti-symmetric index set $I$, the $(n+1-|I|)$-form $\omegas{I}$ is defined by
\begin{gather*}
\omegas{I} = \pm\prod_{a \in \{0,\dots, n\} - I} \omega^a ,
\end{gather*}
with sign so that
\begin{gather*}
\bigg( \prod_{a \in I} \omega^a \bigg) \w \omegas{I} = +\omegas{\varnothing} := \omega^0\w\cdots\w\omega^n .
\end{gather*}
\begin{Theorem}\label{thm:conservation laws linear part}
 For an evolutionary parabolic system $(M,\I)$ with parabolic coframing as in~Pro\-position~$\ref{thm:principal structure equations}$, there are functions $a_i$ and an $(n+1)$-form
 \begin{gather*}
 \Upsilon = \thetas{i,0}\w\omegas{i} - \thetas{\varnothing,0}\w\omegas{0} + a_i\thetas{\varnothing,0}\w\omegas{i}
 \end{gather*}
 so that any element of~$E^{1,n}_1$ is represented by a form
 \begin{gather*}
 \Phi \equiv A\Upsilon - (D_iA) \thetas{\varnothing,0}\w\omegas{i} \mod{\pwf^2},
 \end{gather*}
 where $A$ is a function in $\cMinfty$. The function $A$ satisfies a differential constraint determined by $\dh\Phi \equiv 0$.

 The functions $a_i$ are in $\cinfinity{\M{1}}$ and are determined by the local invariants of~$M_0$. There are functions $W^{I,t}$ ranging over $|I|+t = 2$ so that
 \begin{gather*}
 \dv a_i \equiv V_i^{jkl}\thetas{jkl,0} + W^{I,t}\thetas{Ii,t} \mod{\cartTheta, \pwf^2_2} .
 \end{gather*}
\end{Theorem}

\begin{proof}It suffices to understand the principal weight spectral sequence $\overline{E}_*^{*,k}$ (with $p =1$), which converges to $E_1^{1,*}$.
 From the general theory of characteristic spectral sequences, the zeroth page~$\overline{E}_0$ is isomorphic to the Spencer complex of the tableaux~$K$. This Spencer complex in turn calculates the minimal free resolution of the symbol module associated to $K$. See~\cite{Characteristic_Cohomology_I}, as well as Chapter VIII of~\cite{BCGGG}. The parabolic symbol module has exactly 1 relation, and thus no relations between relations, so the following part of~$\overline{E}_0$ contains the only terms that don't immediately degenerate,
 \[
 \begin{tikzcd}
 \cdots & \overline{E}_0^{n+1,-n+1} & \overline{E}_0^{n+1,-n} & \overline{E}_0^{n+1,-n-1} \\
 \cdots & \overline{E}_0^{n,-n+1} \ar[u, "\dh"] & \overline{E}_0^{n,-n} \ar[u, "\dh"] & 0 \\
 \cdots & \overline{E}_0^{n-1,-n+1} \ar[u, "\dh"] & 0 & 0
 \end{tikzcd}
 \]
Since $\dh$ is a $\cMinfty$-linear vector bundle map between these spaces, it suffices to compute pointwise. Localizing at each point of~$M$ results in the following diagram of vector spaces,
\[
\begin{tikzcd}
\R\{\thetas{ij,0}\w\omegas{\varnothing},\thetas{i,1}\w\omegas{\varnothing},\thetas{\varnothing,2}\w\omegas{\varnothing}\} & \R\{\thetas{i,0}\w\omegas{\varnothing},\thetas{\varnothing,1}\w\omegas{\varnothing}\} & \R\{\thetas{\varnothing,0}\w\omegas{\varnothing}\} \\
\R\{\thetas{i,0}\w\omegas{b},\thetas{\varnothing,1}\w\omegas{b}\} \ar[u, "\dh"] & \R\{\thetas{\varnothing,0}\w\omegas{a}\}\ar[u, "\dh"] & 0 \\
\R\{\thetas{\varnothing,0}\w\omegas{ab}\} \ar[u, "\dh"] & 0 & 0
\end{tikzcd}
\]
Then it is clear that the $\overline{E}_1$ page is
\begin{equation}\label{eq:E_1-bar page spectral sequence}
 \begin{tikzcd}
 0 & 0 & \cMinfty\{\thetas{0,0}\w\omegas{\varnothing}\} \\[-2ex]
 \cMinfty\{\thetas{i,0}\w\omegas{i}\} \ar[ur] & 0 & 0 \\[-2ex]
 0 & 0 & 0
 \end{tikzcd}
\end{equation}
Now there is enough degeneracy to see that the $\overline{E}_2$ page is given by
\[
\begin{tikzcd}
0 & 0 & \cMinfty\{\thetas{0,0}\w\omegas{\varnothing}\} \\[-2ex]
\cMinfty\{\thetas{i,0}\w\omegas{i}\} \ar[urr, "\delta"] & 0 & 0 \\[-2ex]
0 & 0 & 0
\end{tikzcd}
\]
After this page the spectral sequence degenerates, so
\begin{gather*}
E^{1,n}_1 = \ker(\delta) .
\end{gather*}
The operator $\delta$ determines a linear differential operator on functions $A \in \cMinfty$, giving the differential constraint for $A$ to correspond to an element $E^{1,n}_1$.

In~more concrete language, the calculation just done may be unwound as follows. For any element $\Phi\in E^{1,n}_1$, there are functions $A$ and $A_a$ in $\cMinfty$ so that
\begin{gather*}
\Phi \equiv A\,\thetas{i,0}\w\omegas{i} + A_0\,\thetas{\varnothing,0}\w\omegas{0} + A_i\,\thetas{\varnothing,0}\w\omegas{i} \mod{\pwf^2}
\end{gather*}
and $\Phi$ is $\dh$-closed. At highest sub-principal weight, one finds
\begin{gather*}
0 = \dh\Phi \equiv -(A + A_0)\thetas{0,1}\w\omegas{\varnothing} \mod{\thetas{\varnothing,0}, \swf^1_{-n-1}} ,
\end{gather*}
so that $A_0 = -A$. Then, with no restriction on weights,
\begin{gather*}
0 = \dh\Phi \equiv A\dh\left(\thetas{i,0}\w\omegas{i} - \thetas{\varnothing,0}\w\omegas{0}\right) - (D_iA + A_i)\thetas{i,0}\w\omegas{\varnothing} \mod{\thetas{\varnothing,0}, \pwf^2} .
\end{gather*}
Since
\begin{gather*}
\dh\left(\thetas{i,0}\w\omegas{i} - \thetas{\varnothing,0}\w\omegas{0}\right) \equiv 0 \mod{\thetas{\varnothing,0},\thetas{i,0}, \pwf^2} ,
\end{gather*}
there are functions $a_\varnothing, a_i \in \cinfinity{M}$ so that
\begin{gather*}
\dh\left(\thetas{i,0}\w\omegas{i} - \thetas{\varnothing,0}\w\omegas{0}\right) \equiv (a_\varnothing\thetas{\varnothing,0} + a_i\thetas{i,0})\w\omegas{\varnothing} \mod{\pwf^2} ,
\end{gather*}
and it is clear that
\begin{gather*}
A_i = - D_i A + Aa_i .
\end{gather*}
This establishes that
\begin{gather*}
\Phi \equiv A\Upsilon - (D_iA) \thetas{\varnothing,0}\w\omegas{i} \mod{\pwf^2} .
\end{gather*}

Recall the structure equations~\eqref{eq: evolutionary MA structure equations}, which can be restated in the language of weights as
\begin{gather*}
\d\thetas{i,0} \equiv -\beta_i^j\w\thetas{j,0} - \xi_i^j\w\thetas{j,0} - \thetas{ij,0}\w\omega^j - \thetas{i,1}\w\omega^0 \mod{\pwf^2_2},\\
\d\omega^0 \equiv \beta_0^0\w\omega^0 \mod{\omega^0 }, \\
\d\omega^i \equiv - (\delta^i_j\kappa_{\varnothing} - \beta^i_j)\w\omega^j + \xi^i_j\w\omega^j \mod{\omega^0, \pwf^1_0} ,
\end{gather*}
with
\begin{gather*}
\xi_i^j \equiv V_i^{jkl}\thetas{kl,0} \mod{\cartTheta}
\end{gather*}
and
\begin{gather*}
\kappa_\varnothing \equiv \beta_0^0 \equiv \beta_i^j \equiv 0 \mod{\cartTheta, \pwf^1_2} .
\end{gather*}
Using these, one computes
\begin{gather*}
 \d(\thetas{i,0}\w\omegas{i} - \thetanot\w\omegas{0})
 \equiv (a_\varnothing\thetanot + a_i\thetas{i,0})\w\omegas{\varnothing} - 2V_i^{jkl}\thetas{kl,0}\w\thetas{j,0}\w\omegas{i}
 \\ \hphantom{ \d(\thetas{i,0}\w\omegas{i} - \thetanot\w\omegas{0}) \equiv }
 {} + \left(\beta_0^0 -\kappa_{\varnothing} + (n-1)\beta_{\tr}\right)\w\thetas{i,0}\w\omegas{i} \mod{\pwf^2_{2-n}}.
\end{gather*}
Since $\beta_0^0 -\kappa_{\varnothing} + (n-1)\beta_{\tr} \in \pwf^1_2$, there are functions $W^{I,t}$ for $|I|+t = 2$ so that
\begin{gather*}
\beta_0^0 -\kappa_{\varnothing} + (n-1)\beta_{\tr} \equiv -\sum W^{I,t}\thetas{I,t} \mod{\pwf^1_1} .
\end{gather*}
Then
\begin{gather*}
 0 = \d^2(\thetas{i,0}\w\omegas{i} - \thetanot\w\omegas{0})
 \\ \hphantom{ 0 }
 {}\equiv (\d a_i)\w\thetas{i,0}\w\omegas{\varnothing} - 2V_j^{ikl}\thetas{jkl,0}\w\thetas{i}\w\omegas{\varnothing}
 -\!\! \sum\! W^{I,t}\thetas{Ii,t}\w\thetas{i,0}\w\omegas{\varnothing}\mod{\thetanot, \pwf^2_{2-n}}
\end{gather*}
so that, up to an application of Cartan's lemma,
\begin{gather*}
\d a_i \equiv 2V_i^{jkl}\thetas{jkl,0} + \sum W^{I,t}\thetas{Ii,t} \mod{\cartTheta, \pwf^1_2} .
\tag*{\qed}
\end{gather*}
\renewcommand{\qed}{}
\end{proof}

\begin{Remark}\label{rmk: dvUpsilon}
 It will be useful below to note that in the course of the proof, it was shown that
 \begin{align*}
 \d\Upsilon \equiv \dh\Upsilon & - 2V_i^{jkl}(\thetas{kl,0}\w\thetas{j,0} - \thetas{jkl,0}\w\thetanot)\w\omegas{i} \\
 & - W^{I,t}(\thetas{I,t}\w\thetas{i,0} - \thetas{Ii,t}\w\thetanot)\w\omegas{i} \mod{\pwf^2_{2-n}}.
 \end{align*}
\end{Remark}

\section{Conservation laws of strongly parabolic systems}\label{sec:Conservation laws of strong parabolic}
Finally, with the tools developed, it is not difficult to prove the following theorem.
\begin{Theorem}\label{thm:Jacobi potential lives on M_0}
 Let $(M,\I)$ be the infinite prolongation of an evolutionary parabolic system $(M_0,\I_0)$ and
 $\Phi$ a differentiated conservation law on $(M,\I)$. The defining function $A$ of~$\Phi$ factors through $\cinfinity{M_0}$.
\end{Theorem}

\begin{proof}
 By Proposition~\ref{thm:conservation laws linear part}, a conservation law $\Phi$ has a defining function $A \in \cMinfty$ so that the $\I$-linear part of~$\Phi$ is given by
 \begin{gather*}
 \Phi \equiv A\Upsilon - (D_i A)\thetas{\varnothing,0}\w\omegas{i} \mod{\pwf^2} .
 \end{gather*}
 Since $\Phi$ is closed, $\dh\Phi \equiv 0$ and $\dv\Phi$ is defined.
 Calculating directly,
 \begin{gather}\label{eq: dvPhi}
 \dv\Phi \equiv (\d A)\w\Upsilon - \d(D_i A)\w\thetanot\w\omegas{i} \mod{\pwf^2_{2-n},\swf^2_{2-n}} .
 \end{gather}
 Here it is worth noting that only $\I$-quadratic terms need be considered (the $\I$-linear terms cancel each other), and that the $\I$-quadratic terms of~$\d\Upsilon$ were computed in Remark~\ref{rmk: dvUpsilon}.

 Now, using Corollary~\ref{cor:factoring through by weight}, it suffices to demonstrate that
 \begin{gather*}
 \pwt(\d A) \le 2 .
 \end{gather*}
 This is done by first bounding the sub-principal weight of~$\d A$.
 To this end, let
 \begin{gather*}
 w = \wt(\d A) ,
 \end{gather*}
 and $N$ the largest integer for which
 \begin{gather*}
 \dv A \not\equiv 0 \mod{\cartTheta, \pwf^1_{N-1}, \swf^1_{w-1}} .
 \end{gather*}
 Then, for $S = w-N$, there are functions $A_{I,S}$, $B_{J,S+1}$ so that
 \begin{gather*}
 \dv A \equiv \sum_{\substack{|I|+S = N \\ |I|+2S = w}} A_{I,S}\thetas{I,S} + \sum_{\substack{|J|+S+1 = N - 1 \\ |J|+2(S+1) = w}} B_{J,S+1}\thetas{J,S+1} \mod{\cartTheta, \pwf^1_{N-2}, \swf^1_{w-1}} .
 \end{gather*}
 (If $|I| = N-S = 0,1$ then the second sum is vacuous, so let the $B_{J,S+1} = 0$ in that case.)
 From Lemma~\ref{thm:(lem) vertical derivative of directional derivative, sub-principal weights}, if $N \ge 3$
 \begin{gather*}
 \dv(D_i A) \equiv \sum A_{I,S}\thetas{Ii,S} + \sum B_{J,S+1}\thetas{Ji,S+1} \mod{\cartTheta, \pwf^1_{N-1}, \swf^1_{w}} .
 \end{gather*}

 Suppose that either $w \ge 5$ (in which case $N \ge 3$ is automatic) or $w = 4$ and $N \ge 3$, and plug into equation~\eqref{eq: dvPhi},
 \begin{align*}
 \dv\Phi \equiv{}&{ }A_{I,S}\left(\thetas{I,S}\w\thetas{i,0}\w\omegas{i} - \thetas{I,S}\w\thetanot\w\omegas{0} - \thetas{Ii,S}\w\thetas{\varnothing,0}\w\omegas{i}\right) \\
 & + B_{J,S+1}\left(\thetas{J,S+1}\w\thetas{i,0}\w\omegas{i} - \thetas{J,S+1}\w\thetanot\w\omegas{0} - \thetas{Ji,S+1}\w\thetas{\varnothing,0}\w\omegas{i}\right) \\
 & \hspace*{82mm}\mod{\pwf^2_{N-n-1}, \swf^2_{w-n-1}} .
 \end{align*}
 One then computes immediately
 \begin{gather*}
 0 \equiv \dh\dv\Phi \equiv -2A_{I,S}\thetas{I,S+1}\w\thetanot\w\omegas{\varnothing} \mod{\pwf^2_{N-n-1}, \swf^2_{w-n-1}} ,
 \end{gather*}
 so that each $A_{I,S}$ vanishes, contradicting the maximality of~$N$.

 It has just been shown that
 \begin{gather*}
 \dv A \equiv 0 \mod{\cartTheta, \pwf^1_2, \swf^1_3} .
 \end{gather*}
 As such, there are functions $A_{jkl,0}$, $B_{j,1}$, and $B_{\varnothing,2}$ so that
 \begin{gather*}
 \dv A \equiv A_{jkl,0}\thetas{jkl,0} + B_{j,1}\thetas{j,1} + B_{\varnothing,2}\thetas{\varnothing,2} \mod{\cartTheta, \pwf^1_{1}, \swf^1_{2}} .
 \end{gather*}
 The proof of Lemma~\ref{thm:(lem) vertical derivative of directional derivative, sub-principal weights} applies here. Indeed, from
 \begin{gather*}
 0 = \d^2 A \equiv (\dv(D_i A) - A_{jkl,0}\thetas{jkli,0} - B_{j,1}\thetas{ji,1} - B_{\varnothing,2}\thetas{i,2})\w\omega^i \mod{\omega^0, \Lambda^2\cartTheta, \pwf^1_{1}, \swf^1_{2}}
 \end{gather*}
 it is an application of Cartan's lemma to obtain
 \begin{gather*}
 \dv(D_i A) \equiv A_{jkl,0}\thetas{ijkl,0} + B_{j,1}\thetas{ji,1} + B_{\varnothing,2}\thetas{i,2} \mod{\cartTheta, \pwf^1_{2}, \swf^1_{3}} .
 \end{gather*}

 Plugging into equation~\eqref{eq: dvPhi},
 \begin{gather*}
 \dv\Phi \equiv A_{jkl,0}\big(\thetas{jkl,0}\w\thetas{i,0}\w\omegas{i} - \thetas{jkl,0}\w\thetanot\w\omegas{0} - \thetas{ijkl,0}\w\thetas{\varnothing,0}\w\omegas{i}\big)
 \\ \hphantom{\dv\Phi \equiv}
{} + B_{j,1}\big(\thetas{j,1}\w\thetas{i,0}\w\omegas{i} - \thetas{j,1}\w\thetanot\w\omegas{0} - \thetas{ij,1}\w\thetas{\varnothing,0}\w\omegas{i}\big)
\\ \hphantom{\dv\Phi \equiv}
{} + B_{\varnothing,2}\big(\thetas{\varnothing,2}\w\thetas{i,0}\w\omegas{i} - \thetas{\varnothing,2}\w\thetanot\w\omegas{0} - \thetas{i,2}\w\thetas{\varnothing,0}\w\omegas{i}\big)
 \mod{\pwf^2_{2-n}, \swf^2_{2-n}}
 \end{gather*}
 and thus
 \begin{gather*}
 0 \equiv \dh\dv\Phi \equiv -2A_{ijk,0}\thetas{ijk,1}\w\thetanot\w\omegas{\varnothing} \mod{\pwf^2_{2-n}, \swf^2_{2-n}} .
 \end{gather*}
 It follows that each $A_{jkl,0}$ vanishes, which proves the theorem.
\end{proof}

With the theorem just proven, it is not difficult to prove the following theorem.
\begin{Theorem}\label{thm:parabolics with conservation laws are MA}
 Given an evolutionary parabolic system $M_0$ with a non-trivial conservation law~$\Phi$ on its infinite prolongation $M$, in any neighborhood of~$M_0$ where the defining function of~$\Phi$ is not zero, $M_0$ has a Monge--Amp\`ere deprolongation.
\end{Theorem}
\begin{proof}
 Suppose $M_0$ has a non-trivial conservation law
 \begin{gather*}
 \Phi \equiv A\Upsilon - (D_i A)\thetas{\varnothing,0}\w\omegas{i} \mod{\pwf^2} .
 \end{gather*}
 By Theorem~\ref{thm:Jacobi potential lives on M_0}, $A$ is a function on $M_0$. Consequently, in any neighborhood where $A$ is non-zero, there are parabolic cofframings that reduce $A$ to 1. For example, the change of coframing that fixes the $\omega^a$ but transforms the other coframing terms as
 \begin{gather*}
 \thetas{\varnothing,0} \xmapsto{\quad} (1/A)\thetas{\varnothing,0},\qquad\ \,
 \thetas{i,0} \xmapsto{\quad} (1/A)\thetas{i,0}, \qquad
 \thetas{\varnothing,1} \xmapsto{\quad} (1/A)\thetas{\varnothing,1} ,
 \\
 \pis{ij,0} \xmapsto{\quad} (1/A)\pis{ij,0}, \qquad\ \,
 \pis{i,1} \xmapsto{\quad} (1/A)\pis{i,1}, \qquad
 \pis{\varnothing,2} \xmapsto{\quad} (1/A)\pis{\varnothing,2}
 \end{gather*}
 preserves the evolutionary parabolic structure of~$(M_0,\I_0)$, and after this reduction,
 \begin{gather*}
 \Phi \equiv \Upsilon \mod{\pwf^2} .
 \end{gather*}

 By Remark~\ref{rmk: dvUpsilon},
 \begin{gather*}
 \dv\Phi \equiv - 2V_i^{jkl}(\thetas{kl,0}\w\thetas{j,0} - \thetas{jkl,0}\w\thetanot)\w\omegas{i}
 \\ \hphantom{\dv\Phi \equiv}
{} - W^{I,t}(\thetas{I,t}\w\thetas{i,0} - \thetas{Ii,t}\w\thetanot)\w\omegas{i}
 \mod{\pwf^2_{2-n}} .
 \end{gather*}
 Then, recalling that $V_i^{jkl}$ has values in $\Sym_0^2 (\Sym_0^2 \R^n)$,
 \begin{gather*}
 0 = \dh\dv\Phi \equiv -2V_i^{jkl}\thetas{ijkl,0}\w\thetanot\w\omegas{\varnothing} \mod{\pwf^2_{2-n}} ,
 \end{gather*}
 and the $\Sym_0^4 \R^n$ component of the Monge--Amp\`ere invariant vanishes identically.
\end{proof}

\subsection*{Acknowledgements}
This material is based upon work supported by the National Science Foundation under Grants No.~DGE-1106400 and 74341.2010, as well as the Australian Research Council, Discovery Program DP190102360.

\pdfbookmark[1]{References}{ref}
\LastPageEnding

\end{document}